\documentclass[12pt,a4paper,oneside,reqno,notitlepage]{amsart}
\usepackage{amssymb}
\usepackage{color}
\textwidth
150mm

\theoremstyle{plain}

\numberwithin{equation}{section}

\begin{document}

\baselineskip 8mm
\parindent 9mm

\title[]
{Fractional order semilinear Volterra integrodifferential equations in Banach spaces}

\author{Kexue  Li}

\address{School of Mathematics and Statistics, Xi'an Jiaotong University, Xi'an 710049, China}

\email{kexueli@gmail.com}


\begin{abstract}
In this paper, sufficient conditions are established for the existence results of fractional order semilinear Volterra integrodifferential equations in Banach spaces. The results are obtained by using the theory of fractional cosine families and fractional powers of operators.
\end{abstract}
\maketitle

\section{\textbf{Introduction}}
The integrodifferential equations in Banach spaces have attracted much interest.
Pr$\ddot{\mbox{u}}$ss \cite{Pruss} considered the solvability behavior on the real line of linear integrodifferential equations in a general Banach space and gave several applications to integral partial differential equations. Grimmer \cite{Grimmer} established general conditions to ensure the existence of a resolvent operator for an integrodifferential equation in a Banach space. Fitzgibbon \cite{WEF} studied the existence, continuation, and behavior of solutions to an abstract semilinear Volterra integrodifferential equation.
Keyantuo and Lizama \cite{Valentin} characterized existence and uniqueness of solutions for a linear integro-differential equation
 in  Holder spaces. Londen \cite{Londen} proved an existence result on a nonlinear Volterra integrodifferential equation in a real reflexive Banach space by using  the theory of maximal monotone operators.
 Pr$\ddot{\mbox{u}}$ss \cite{J} studied linear Volterra integrodifferential equations in a Banach space in case the main part of the equation generates an analytic $C_{0}$-semigroup.  Travis and Webb \cite{Travis} studied the existence of solutions to semilinear second order Volterra integrodifferential equations in Banach spaces by using the theory of strongly continuous cosine families.  Mainini and Mola \cite{MM} considered in an abstract setting, an instance of the
Coleman-Gurtin model for heat conduction with memory. Engler \cite{Engler} constructed global weak solution  of  scalar  second-order quasilinear hyperbolic integrodifferential equations  with singular kernels.
 Pr$\ddot{\mbox{u}}$ss \cite{Jan} studied the existence, positivity, regularity, compactness and integrability of the resolvent for a class of Volterra equations of scalar type.
 Hern$\acute{\mbox{a}}$ndez \cite{EH} studied the existence of strict and classical solutions for a class of abstract non-autonomous Volterra integrodifferential equations in Banach spaces. Lang and Chang \cite{Lang} investigated the local existence and uniqueness of solutions to integrodifferential equations with infinite delay. Jawahdou \cite{Adel} studied the existence of mild solutions for initial value problems for semilinear Volterra integrodifferential equations in a Banach space.

  In recent years, fractional differential equations have received increasing attention due to its applications in physics, chemistry, materials, engineering, biology, finance, we refer to \cite{Pod}, \cite{Main}, \cite{R}. Fractional order derivatives  have the memory property and can describe many phenomena that integer order derivatives can't characterize.

Consider the following fractional semilinear differential equation
\begin{equation}\label{cauchy}
 \left\{\begin{aligned} &^CD_t^\alpha u(t)=Au(t), \, t>0;\\
&u(0)=x, u^{(k)}(0)=0, \ k=1,2,\cdots, m-1, \end{aligned}\right.
\end{equation}
where $\alpha>0$, $m$ is the smallest integer greater than or equal
to $\alpha$, $^CD_t^\alpha$ is the $\alpha$-order Caputo fractional derivative operator, $A: D(A)\subset X\rightarrow X$ is a closed densely defined linear operator on a Banach space $X$.

Bazhlekova \cite{bazh} introduced the notion of solution operator for (\ref{cauchy}) as follows.\\
\textbf{Definition 1.1.}
A family $\{C_{\alpha}(t)\}_{t\geq0}\subset \mathcal{B}(X)$ is called a solution operator for (\ref{cauchy}) if the following conditions are satisfied:\\
(a) $C_{\alpha}(t)$ is strongly continuous for $t\geq 0$ and $C_{\alpha}(0)=I$ (the identity operator on $X$);\\
(b) $C_{\alpha}(t)D(A)\subset D(A)$ and $AC_{\alpha}(t)\xi=C_{\alpha}(t)A\xi$ for all $\xi\in D(A)$, $t\geq 0$;\\
(c) $C_{\alpha}(t)\xi$ is a solution of $x(t)=\xi+\int_{0}^{t}g_{\alpha}(t-s)Ax(s)ds$ for all $\xi\in D(A)$, $t\geq 0$,\\
we refer to the equality (2.3) about the definition of $g_{\alpha}(t)$.

 $A$ is called the infinitesimal generator of $C_{\alpha}(t)$.  Note that in some literature the solution operator  also is called fractional resolvent family or fractional resolvent operator function, we refer to \cite{Chenchuang}, \cite{Miao}. As a matter of fact, the solution operator $C_{2}(t)$ is a cosine family, in this paper, for $\alpha\in(1,2]$, the solution operator $C_{\alpha}(t)$ is called strongly continuous $\alpha$-order fractional cosine family, or $\alpha$-order cosine family, for short.

Chen and Li \cite{Chenchuang} developed a purely algebraic notion, named
$\alpha$-resolvent operator function: A family
$\{S_\alpha(t)\}_{t\geq 0}$ of bounded linear operators of a Banach space $X$ is
called an $\alpha$-resolvent operator function if the following
conditions are satisfied:

(a) $S_\alpha(t)$ is strongly continuous for $t\geq 0$ and
$S_\alpha(0)=I\ (\mbox{the identity operator})$,

(b) $S_\alpha(t)S_\alpha(s)=S_\alpha(s)S_\alpha(t)$ for all $t,
s\geq 0$, and

(c) there holds for all $t, s\geq 0$ that
\begin{equation*}
S_\alpha(s)J^\alpha_tS_\alpha(t)-J^\alpha_sS_\alpha(s)S_\alpha(t)=J^\alpha_tS_\alpha(t)-J^\alpha_sS_\alpha(s),
\end{equation*}
where $J^\alpha_t$ is the $\alpha$-order Riemann-Liouville fractional
integral operator.

It has been proved in \cite{Chenchuang} that a family
$\{S_\alpha(t)\}_{t\geq 0}$ is an $\alpha$-resolvent operator
function if and only if it is a solution operator (or an
$\alpha$-times resolvent family, see \cite{Miao}) for the problem (\ref{cauchy}).

Peng and Li \cite{Peng} developed a novel operator theory for the
problem (\ref{cauchy}) with the order $\alpha\in(0,1)$.\\
\textbf{Definition
1.1.} (\cite{Peng}). Let $0<\alpha<1$. A one-parameter family $\{T_{\alpha}(t)\}_{t\geq 0}$ of bounded
linear operators of $X$ is called strongly continuous fractional
semigroup of order $\alpha$ (or $\alpha$-order fractional semigroup,
for short) if it possesses the following two properties:

$(I)$ for every $x\in X$, the mapping $t\mapsto T(t)x$ is continuous
over $[0,\infty)$;

$(II)$ $T_{\alpha}(0)=I$, and for all
$t,s\geq 0$,
\begin{equation}\label{eq01}
\begin{aligned}
&\int^{t+s}_0
\frac{T_{\alpha}(\tau)d\tau}{(t+s-\tau)^{\alpha}}-\int^{t}_{0}\frac{T_{\alpha}(\tau)d\tau}{(t+s-\tau)^{\alpha}}-\int^{s}_{0}\frac{T_{\alpha}(\tau)d\tau}{(t+s-\tau)^{\alpha}}\\
=&\alpha\int^t_0\int^s_0\frac{T_{\alpha}(\tau_1)T_{\alpha}(\tau_2)}{(t+s-\tau_1-\tau_2)^{1+\alpha}}d\tau_1d\tau_2,
\end{aligned}
\end{equation}
where the integrals are understood in the sense of strong
operator topology.

For $\alpha\in(0,1)$, it is proved that a family
of bounded linear operator is a solution operator for (\ref{cauchy}) if and
only if it is a fractional semigroup. Moreover, it is shown that the problem (\ref{cauchy}) is well-posed if and only if its coefficient operator generates an $\alpha$-order semigroup.

Keyantuo \cite{Keyantuo} investigated a general framework for connections
between ordinary non-homogeneous equations in Banach spaces and fractional Cauchy problems. When the underlying operator generates a strongly continuous semigroup, using a subordination argument, the fractional evolution equation is well-posed.

In this paper we are concerned with the  fractional order semilinear Volterra integrodifferential equation
\begin{equation}\label{1}\left\{
\begin{array}{ll}
^{C}D_{t}^{\alpha}u(t)=Au(t)+\int_{0}^{t}h(t,s,u(s))ds+f(t), \quad t\in R_{+},\vspace{2mm}\\
u(0)=x, \quad u'(0)=y,
\end{array} \right.
\end{equation}
where $R_{+}=[0,\infty)$, $\alpha\in(1,2]$, $^{C}D_{t}^{\alpha}$ is the $\alpha$-order Caputo fractional derivative operator, $A$ is the infinitesimal generator of a strongly continuous  fractional cosine family $\{C_{\alpha}(t)\}_{t\geq 0}$ on a Banach space $X$, $h$ is a nonlinear unbounded operator from $R_{+}\times R_{+}\times X$ to $X$, $f$ is a function from $R_{+}$ to $X$ and $x,y\in X$.

This paper is organized as follows. In Section 2, we give the basic notations and preliminary facts. In Section 3, we give the sufficient conditions  for the existence of the equation (\ref{1}).  At last, an example is presented to illustrate the main results.
\section{Preliminaries}
Let $X$ be a Banach space with norm $\|\cdot\|$. By $\mathcal{B}(X)$ we denote the space of all bounded linear operators on $X$. Let $1\leq p<\infty$. By $L^{p}([0,T]; X)$ we denote the space of
$X$-valued Bochner integrable functions $f:[0,T]\rightarrow X$ with
the norm
\begin{equation}
\|f\|_{L^{p}([0,T];X)}=(\int_{0}^{T}\|f(t)\|^{p}dt)^{1/p}
\end{equation}
By $C([0,T];X)$, resp. $C^{1}([0,T];X)$, we denote the spaces of functions $f:[0,T]\rightarrow X$, which are continuous, resp. $1$-times continuously differentiable. $C([0,T];X)$ and $C^{1}([0,T];X)$ are Banach spaces endowed with the norms
\begin{equation}
\|f\|_{C}=\sup_{t\in [0,T]}\|f(t)\|_{X}, \quad ||f\|_{C^{1}}=\sup_{t\in [0,T]}\sum_{k=0}^{1}\|f^{(k)}(t)\|_{X}.
\end{equation}
Let $I$ be the identity operator on $X$. If $A$ is a linear operator on $X$, then $R(\lambda,A)=(\lambda I-A)^{-1}$ denotes the resolvent operator of $A$.
For the sake of simplicity, we use the notation for $\alpha>0$:
\begin{equation}
g_{\alpha}(t)=\frac{t^{\alpha-1}}{\Gamma(\alpha)},\quad t>0,
\end{equation}
where $\Gamma(\alpha)$ is the Gamma function. If $\alpha=0$, we set
$g_{0}(t)=\delta(t)$, the delta distribution.\\
\textbf{Definition 2.1}. The Riemann-Liouville fractional integral
 of order $\alpha\geq0$ is defined by
\begin{equation}
J_{t}^{\alpha}u(t)=\int_{0}^{t}g_{\alpha}(t-s)u(s)ds,
\end{equation}
where $u(t)\in L^{1}([0,T];X)$.

The set of the Riemann-Liouville fractional integral operators $\{J_{t}^{\alpha}\}_{\alpha\geq 0}$ is a semigroup, i.e., $J_{t}^{\alpha}J_{t}^{\beta}=J_{t}^{\alpha+\beta}$ for all $\alpha,\beta\geq 0$. \\
\textbf{Definition 2.2}. The Riemann-Liouville fractional derivative
of order $\alpha\in(1,2]$ is defined by
\begin{equation} \label{0}
 D_{t}^{\alpha}u(t)=\frac{d^{2}}{dt^{2}}J_{t}^{2-\alpha}u(t),
\end{equation}
where $u(t)\in L^{1}([0,T];X)$, \ $D_{t}^{\alpha}u(t)\in L^{1}([0,T];X)$.\\
\textbf{Definition 2.3}. The Caputo fractional derivative
of order $\alpha\in(1,2]$ is defined by
\begin{equation} \label{0}
 ^{C}D_{t}^{\alpha}u(t)=D_{t}^{\alpha}(u(t)-u(0)-u'(0)t),
\end{equation}
where $u(t)\in L^{1}([0,T];X)\cap C^{1}([0,T];X)$, \ $D_{t}^{\alpha}u(t)\in L^{1}([0,T];X)$.\\
The Laplace transform for the Riemann-Liouville fractional
integral is given by
\begin{equation}\label{well}
L\{J_{t}^{\alpha}u(t)\}=\frac{1}{\lambda^{\alpha}}\widehat{u}(\lambda),
\end{equation}
where $\widehat{u}(\lambda)$ is the Laplace of $u$ given by
\begin{equation}\label{well2}
\widehat{u}(\lambda)=\int_{0}^{\infty}e^{-\lambda t}u(t)dt,\ \mathrm{Re}
\lambda>\omega.
\end{equation}
The Laplace transform for Caputo derivative is given by
\begin{equation}\label{mild solution}
L\{^{C}D_{t}^{\alpha}u(t)\}=\lambda^{\alpha}\widehat{u}(\lambda)-u(0)\lambda^{\alpha-1}-u'(0)\lambda^{\alpha-2}.
\end{equation}
\textbf{Definition 2.4.} The fractional sine family $S_{\alpha}: R_{+}\rightarrow \mathcal{B}(X)$ associated with $C_{\alpha}$ is defined by
\begin{equation}\label{sine}
S_{\alpha}(t)=\int_{0}^{t}C_{\alpha}(s)ds.
\end{equation}
\textbf{Remark 2.5.} For $x\in X$, Define $S'(0)x=\frac{dS_{\alpha}(t)x}{dt}|_{t=0}$. From Definition 2.5 and Definition 1.1, it is clear that $S'(0)=I\ (\mbox{the identity operator on X})$.\\
\textbf{Definition 2.6.} The fractional Riemann-Liouville family $P_{\alpha}: R_{+}\rightarrow \mathcal{B}(X)$  associated with $C_{\alpha}$ is defined by
\begin{equation}\label{RL}
P_{\alpha}(t)=J_{t}^{\alpha-1}C_{\alpha}(t).
\end{equation}
\textbf{Definition 2.7.} The $\alpha$-order cosine family $C_{\alpha}(t)$ is called exponentially bounded if there are constants $M\geq 1$ and $\omega\geq 0$ such that
\begin{equation}\label{def2.7}
\|C_{\alpha}(t)\|\leq Me^{\omega t},\ t\geq 0.
\end{equation}
An operator $A$ is said to belong to $\mathcal{C^{\alpha}}(M,\omega)$, if the problem (\ref{cauchy}) has an $\alpha$-order cosine family  $C_{\alpha}(t)$ satisfying (\ref{def2.7}).
\section{Existence of solutions}
For $\alpha\in (1,2)$, we assume $A\in \mathcal{C^{\alpha}}(M,\omega)$ and let $C_{\alpha}(t)$ be the corresponding $\alpha$-order cosine family.
We have (see \cite{bazh}, (2.5) and (2.6))
\begin{equation}\label{spectrum}
\{\lambda^{\alpha}: \mbox{Re}\lambda>\omega\}\subset \rho(A),
\end{equation}
and
\begin{equation}\label{resolvent}
\lambda^{\alpha-1}R(\lambda^{\alpha},A)\xi=\int_{0}^{\infty}e^{-\lambda t}C_{\alpha}(t)\xi dt, \ \mathrm{Re}\lambda>\omega,\ \xi\in X.
\end{equation}
By (\ref{RL}), (\ref{resolvent}), we have
\begin{equation}\label{rho}
R(\lambda^{\alpha},A)\xi=\int_{0}^{\infty}e^{-\lambda t}P_{\alpha}(t)\xi dt,\ \mathrm{Re}\lambda>\omega,\ \xi\in X.
\end{equation}
For a fractional cosine family $C_{\alpha}(t)$, we define $E=\{x\in X: \ C_{\alpha}(t)x$ is once continuously
differentiable on $R_{+}\}$. By the identity $\lambda^{\alpha}R(\lambda^{\alpha},A)-I=AR(\lambda^{\alpha},A)$,  (\ref{resolvent}) and (\ref{rho}), we have that $P_{\alpha}(t)E\subset D(A), \ t\in R_{+}$, and
\begin{equation}\label{key}
\frac{d}{dt}C_{\alpha}(t)x=AP_{\alpha}(t)x, \ x\in E, \ t\in R_{+}.
\end{equation}
By (\ref{sine}), (\ref{resolvent}), we have
\begin{equation}\label{lambda}
\lambda^{\alpha-2}R(\lambda^{\alpha},A)\xi=\int_{0}^{\infty}e^{-\lambda t}S_{\alpha}(t)\xi dt,\quad \mathrm{Re}\lambda>\omega, \ \xi\in X.
\end{equation}
\textbf{Lemma 3.1.} Let $A$ be the infinitesimal generator of an
$\alpha$-order cosine family $C_{\alpha}(t)$, and $S_{\alpha}(t)$ is the corresponding $\alpha$-order sine family.
Then \\
(a) For all $x\in D(A)$ and $t\geq0$,
 \begin{equation*}
S_{\alpha}(t)x\in D(A)\ \ \mbox{and} \ \ AS_{\alpha}(t)x=S_{\alpha}(t)Ax.
 \end{equation*}
(b) For all $x\in D(A)$ and $t\geq0$,
 \begin{equation*}
 S_{\alpha}(t)x=tx+J_{t}^{\alpha}S_{\alpha}(t)Ax.
 \end{equation*}
\textbf{Proof.} (a) Fix some  $\mu^{\alpha}\in \rho(A)$, for $\lambda>\max\{\omega,0\}$ and $x\in X$,
 \begin{align*}
\int_{0}^{\infty}e^{-\lambda t}S_{\alpha}(t)\mu^{\alpha-2}R(\mu^{\alpha},A)xdt &=\lambda^{\alpha-2}R(\lambda^\alpha,A)\mu^{\alpha-2}R(\mu^\alpha,A)x\\
 & = \mu^{\alpha-2}R(\mu^\alpha,A)\lambda^{\alpha-2}R(\lambda^\alpha,A)x\\
 &=\int_{0}^{\infty}e^{-\lambda t}\mu^{\alpha-2}R(\mu^{\alpha},A)S_{\alpha}(t)xdt.
\end{align*}
From the uniqueness theorem of the Laplace transform, it follows that
$R(\mu^{\alpha},A)S_{\alpha}(t)=S_{\alpha}(t)R(\mu^{\alpha},A)$. This implies (a).\\
(b) For $x\in D(A)$, $\lambda>\omega\geq0$,
\begin{align*}
\int_{0}^{\infty}\lambda^{2}e^{-\lambda t}txdt&=x\nonumber\\
&=\lambda^{\alpha}R(\lambda^{\alpha},A)x-R(\lambda^{\alpha},A)Ax\nonumber\\
&=\int_{0}^{\infty}\lambda^{2}e^{-\lambda t}S_{\alpha}(t)xdt-\int_{0}^{\infty}\lambda^{2}e^{-\lambda t}J_{t}^{\alpha}S_{\alpha}(t)Axdt.
\end{align*}
Hence, (b) follows from the uniqueness theorem of Laplace transforms. \ \ \ $\Box$

Since $A\in \mathcal{C^{\alpha}}(M,\omega)$ for $\alpha\in (1,2)$, then from Theorem 3.3 in \cite{bazh}, it follows that $A$ generates an analytic semigroup $T(t)$ of angle $(\alpha-1)\pi/2$. We suppose that $0\in \rho(A)$, then for $\beta\in (0,1)$, we can define the fractional powers operator $(-A)^{-\beta}$ as follows:
\begin{equation*}
(-A)^{-\beta}=\frac{\sin \pi \beta}{\pi}\int_{0}^{\infty}\tau^{-\beta}(\tau I-A)^{-1}d\tau.
\end{equation*}
\textbf{Definition 3.2.}
Let $A$ be the infinitesimal generator of an analytic semigroup $T(t)$. For every $\beta>0$ we define
\begin{equation*}
(-A)^{\beta}=((-A)^{-\beta})^{-1}.
\end{equation*}
For $\beta=0$, $(-A)^{\beta}=I$.

We collect some basic properties of fractional powers $(-A)^{\beta}$ in the following.\\
\textbf{Lemma 3.3.} (\cite{Pazy})
Assume $(-A)^{\beta}$ is defined by Definition 3.1, then\\
(a) $(-A)^{\beta}$ is a closed operator with domain $D((-A)^{\beta})=R((-A)^{-\beta})$ (the range of $(-A)^{-\beta})$. \\
(b) For $\beta\geq \gamma>0$, $D((-A)^{\beta})\subset D((-A)^{\gamma})$. \\
(c) $D((-A)^{\beta})$ is dense in $X$ for every $\beta\geq 0$. \\
(d) If $\beta,\gamma$ are real then
\begin{equation*}
(-A)^{\beta+\gamma}x=(-A)^{\beta}(-A)^{\gamma}x
\end{equation*}
for every $x\in D((-A)^{\eta})$ where $\eta=\max (\beta,\gamma,\beta+\gamma)$.

By (c),(d) of Lemma 3.3, we see that for $\beta\in (0,1)$,
\begin{equation}\label{beta}
(-A)^{\beta}=(-A)^{\beta-1}(-A).
\end{equation}

We note that $D((-A)^{\beta})$ is a Banach space equipped with the norm $\|x\|_{\beta}=\|(-A)^{\beta}x\|$, $x\in D((-A)^{\beta})$.  By $X_{\beta}$ we denote this Banach space. \\
\textbf{Lemma 3.4.}
Let $A$ be the infinitesimal generator of an $\alpha$-order cosine family $C_{\alpha}(t)$ on $X$. By $P_{\alpha}(t)$ we denote the corresponding Riemann-Liouville family. If $k: R_{+}\rightarrow X$ is continuously differentiable and $v(t)=\int_{0}^{t}P_{\alpha}(t-s)k(s)ds$, then
$v(t)\in D(A)$ for $t\geq 0$, and
\begin{align}\label{lemma}
Av(t)=\int_{0}^{t}C_{\alpha}(t-S)k'(\tau)ds+C_{\alpha}(t)k(0)-k(t).
\end{align}
\textbf{Proof.} Since $k: R_{+}\rightarrow X$ is continuously differentiable, we have
\begin{align}\label{vt}
v(t)&=\int_{0}^{t}P_{\alpha}(t-s)k(s)ds\nonumber\\
&=\int_{0}^{t}P_{\alpha}(t-s)(\int_{0}^{s}k'(\tau)d\tau+k(0))ds\nonumber\\
&=\int_{0}^{t}\int_{0}^{t-\tau}P_{\alpha}(s)k'(\tau)dsd\tau+\int_{0}^{t}P_{\alpha}(t-s)k(0)ds.
\end{align}
From (\ref{RL}), (b) of Proposition 3.3 in \cite{Chenchuang}, it follows that for all
$x\in X$, $t\geq 0$,
$\int_{0}^{t}P_{\alpha}(s)xds\in D(A)$ and
\begin{align*}
A\int_{0}^{t}P_{\alpha}(s)xds=C_{\alpha}(t)x-x.
\end{align*}
Then $v(t)\in D(A)$,
\begin{align}\label{avt}
Av(t)&=\int_{0}^{t}(C_{\alpha}(t-\tau)k'(\tau)-k'(\tau))d\tau+C_{\alpha}(t)k(0)-k(0)\nonumber\\
&=\int_{0}^{t}C_{\alpha}(t-s)k'(s)ds+C_{\alpha}(t)k(0)-k(t). \ \ \ \Box
\end{align}
\textbf{Lemma 3.5.} Let $A$ be the infinitesimal  generator of  an $\alpha$-order cosine family $C_{\alpha}(t)$ on $X$. If $f: R_{+}\rightarrow X$ is continuously differentiable, $x,y\in D(A)$, and let $\varphi(t)=C_{\alpha}(t)x+S_{\alpha}(t)y+\int_{0}^{t}P_{\alpha}(t-s)f(s)ds$, $t\in[0,T]$, then $\varphi(t)\in D(A)$ and $\varphi$ satisfies
\begin{equation*}
 \ \left\{\begin{aligned} &^{C}D_{t}^{\alpha}\varphi(t)=A\varphi(t)+f(t), \ t\in R_{+},\\
&\varphi(0)=x, \ \varphi'(0)=y.
\end{aligned}\right.
\end{equation*}
\textbf{Proof.} From (\ref{key}) and Lemma 3.3, it follows that $\varphi(t)\in D(A)$. It is clear that $\varphi(0)=x$. Since $f: R_{+}\rightarrow X$ is continuously differentiable, it is easy to show that $\varphi'(0)=y$. By (\ref{0}), Remark 2.5 and Lemma 3.1, we have
\begin{align*}
^{C}D_{t}^{\alpha}\varphi(t)&=^{C}D_{t}^{\alpha}C_{\alpha}(t)x+^{C}D_{t}^{\alpha}S_{\alpha}(t)y+^{C}D_{t}^{\alpha}(\int_{0}^{t}P_{\alpha}(t-s)f(s)ds)\\
&=AC_{\alpha}(t)x+D_{t}^{\alpha}(S_{\alpha}(t)y-S_{\alpha}(0)y-tS_{\alpha}'(0)y)+D_{t}^{\alpha}(\int_{0}^{t}P_{\alpha}(t-s)f(s)ds)\\
&=AC_{\alpha}(t)x+D_{t}^{\alpha}(S_{\alpha}(t)y-ty)+\frac{d^{2}}{dt^{2}}J_{t}^{2-\alpha}(P_{\alpha}(t)\ast f(t))\\
&=AC_{\alpha}(t)x+D_{t}^{\alpha}J_{t}^{\alpha}S_{\alpha}(t)Ay+\frac{d^{2}}{dt^{2}}(g_{2-\alpha}(t)\ast g_{\alpha-1}(t)\ast C_{\alpha}(t)\ast f(t))\\
&=AC_{\alpha}(t)x+S_{\alpha}(t)Ay+\frac{d^{2}}{dt^{2}}(1\ast C_{\alpha}(t)\ast f(t))\\
&=AC_{\alpha}(t)x+AS_{\alpha}(t)y+\frac{d}{dt}(C_{\alpha}(t)\ast f(t)).
\end{align*}
By Lemma 3.4, we have
\begin{align*}
\frac{d}{dt}(C_{\alpha}(t)\ast f(t))=A\int_{0}^{t}P_{\alpha}(t-s)f(s)ds+f(t).
\end{align*}
Therefore, the proof is complete. \ \ \ $\Box$

We make the following assumptions on the functions $h$ and $f$:\\
($A_{1}$) $h:R_{+}\times R_{+} \times  D \rightarrow X$ is continuous, where $D$ is an open subset of $X_{\beta}$, $\beta\in [0,1)$.\\
($A_{2}$) $h_{1}: R_{+}\times R_{+} \times  D \rightarrow X$ is continuous, where $h_{1}$ denotes the derivative of $h$ with respect to its first variable.\\
($A_{3}$) $f: R_{+}\rightarrow X$ is continuously differentiable.\\
\textbf{Theorem 3.6.} Let $\alpha\in (1,2)$. Assume that $A\in \mathcal{C^{\alpha}}(M,\omega)$ and let $C_{\alpha}(t)$, $S_{\alpha}(t)$ and  $P_{\alpha}(t)$ denote the corresponding $\alpha$-order cosine family, $\alpha$-order sine family and $\alpha$-order Riemann-Liouville family, respectively.  Assume that $A^{-1}$ is compact. Let $x\in D$, $\beta\in (0,1)$ and let $(-A)^{\beta-1}y\in E$. If ($A_{1}$), ($A_{2}$) and ($A_{3}$) are satisfied, then there exists $T>0$ and a continuous function $u:[0,T]\rightarrow X_{\beta}$ such that
\begin{equation}\label{ut}
u(t)=C_{\alpha}(t)x+S_{\alpha}(t)y+\int_{0}^{t}P_{\alpha}(t-s)\int_{0}^{s}h(s,r,u(r))drds+\int_{0}^{t}P_{\alpha}(t-s)f(s)ds, \ t\in[0,T].
\end{equation}
If, in addition, $x\in D(A)$ and $y\in E$, then the Caputo derivative $^{C}D_{t}^{\alpha}u$ of the solution $u$ of (\ref{ut}) is continuous, $u\in D(A)$, and $u$ satisfies
\begin{equation}\label{Caputo}
 \ \left\{\begin{aligned} &^{C}D_{t}^{\alpha}u(t)=Au(t)+\int_{0}^{t}h(t,s,u(s))ds+f(t), \ t\in [0,T],\\
&u(0)=x, \ u'(0)=y.
\end{aligned}\right.
\end{equation}
\textbf{Proof.} For $\delta>0$, let $N_{\delta}(x)=\{x_{1}\in X_{\beta}: \|x-x_{1}\|_{\beta}< \delta\}$. Let $\varphi(t)=C_{\alpha}(t)x+S_{\alpha}(t)y+\int_{0}^{t}P_{\alpha}(t-s)f(s)ds$. We can choose $\delta>0$ and $T>0$ such that $N_{\delta}(x)\subset D$ and
for $r,s\in [0,T]$ and $x_{1}\in N_{\delta}(x)$,
\begin{equation}\label{c1}
\|h(r,s,x_{1})\|\leq 1+M(x,T), \ \|h_{1}(r,s,x_{1})\|\leq1+N(x,T);
\end{equation}
for $t\in [0,T]$,
\begin{equation}\label{c2}
\|\varphi(t)-x\|_{\beta}< \delta/2;
\end{equation}
for $t \in[0,T]$ and $x_{1},x_{2},x_{3}\in N_{\delta}(x)$,
\begin{equation}\label{c3}
\|(-A)^{\beta-1}\left(-\int_{0}^{t}C_{\alpha}(t-s)(h(s,s,x_{1})+\int_{0}^{s}h_{1}(s,r,x_{2})dr)ds+\int_{0}^{t}h(t,s,x_{3})ds\right)\|<\frac{\delta}{2},
\end{equation}
where $M(x,T)=\sup_{r,s\in [0,T]}\|h(r,s,x)\|$ and $N(x,T)=\sup_{r,s\in [0,T]}\|h_{1}(r,s,x)\|$.\\
In fact, since $h:R_{+}\times R_{+} \times  D \rightarrow X$ is continuous, given $\varepsilon>0$, there exists $\delta>0$ such that for $x_{1} \in N_{\delta}(x)$ and $r,s\in [0,T]$, we have
\begin{equation*}
\|h(s,r,x_{1})-h(s,r,x)\|<\varepsilon.
\end{equation*}
Letting  $\varepsilon\in (0,1)$,  we obtain $\|h(r,s,x_{1})\|\leq 1+M(x,T)$. Similiarly, $\|h_{1}(r,s,x_{1})\|\leq 1+N(x,T)$. \\
It is easy to show that $(-A)^{\beta}C_{\alpha}(t)x=C_{\alpha}(t)(-A)^{\beta}x$ for $x\in X_{\beta}$.
Note that $t\in [0,T]$ and $C_{\alpha}(t)$ is strongly continuous for $t\geq 0$, then
\begin{align}\label{change}
\|(-A)^{\beta}(C_{\alpha}(t)x-x)\|&=\|(C_{\alpha}(t)-I)(-A)^{\beta}x\|\nonumber\\
& \leq \|C_{\alpha}(t)-I\|\|x\|_{\beta}\nonumber\\
&=C(\alpha,T)\|x\|_{\beta},
\end{align}
where $C(\alpha,T)=\sup_{t\in [0,T]}\|C_{\alpha}(t)-I\|$. \\
Since $\beta\in (0,1)$, then there exists a positive constant $M_{0}>0$ such that
\begin{align*}
\|(-A)^{\beta-1}\|\leq M_{0},
\end{align*}
see Lemma 6.3 in \cite{Pazy}. \\
Since $(-A)^{\beta-1}y\in E$, we have
\begin{align}\label{ca}
\|(-A)^{\beta}S_{\alpha}(t)y\|&=\|(-A)^{\beta-1}(-A)S_{\alpha}(t)y\|\nonumber\\
&=\|(-A)^{\beta-1}\frac{d}{dt}C_{\alpha}(t)y\|\nonumber\\
&\leq \|(-A)^{\beta-1}\|\|\frac{d}{dt}C_{\alpha}(t)y\|\nonumber\\
&\leq M_{0}M'(\alpha,T,y),
\end{align}
where $M'(\alpha,T,y)=\sup_{t\in [0,T]}\frac{d}{dt}C_{\alpha}(t)y$. \\
By Lemma 3.4, we have
\begin{align}\label{ap}
\|(-A)^{\beta}\int_{0}^{t}P_{\alpha}(t-s)f(s)ds\|&=\|(-A)^{\beta-1}(-A)\int_{0}^{t}P_{\alpha}(t-s)f(s)ds\|\nonumber\\
&\leq \|(-A)^{\beta-1}\|\|(-A)\int_{0}^{t}P_{\alpha}(t-s)f(s)ds\|\nonumber\\
&=\|(-A)^{\beta-1}\|\|\int_{0}^{t}C_{\alpha}(t-s)f'(s)ds+C_{\alpha}(t)f(0)-f(t)\|\nonumber\\
&=\|(-A)^{\beta-1}\|(M'_{T}Me^{\omega T}T+M_{T})\nonumber\\
&\leq M_{0}(M'_{T}Me^{\omega T}T+M_{T}).
\end{align}
where $M'_{T}=\sup_{s\in [0,T]}\|f'(s)\|$, $M_{T}=\sup_{s\in [0,T]}\|C_{\alpha}(t)f(0)-f(t)\|$.
Since
\begin{align}\label{varphi}
\|\varphi(t)-x\|_{\beta}&=\|(-A)^{\beta}(\varphi(t)-x)\|\nonumber\\
&\leq \|(-A)^{\beta}(C_{\alpha}(t)x-x)\|+\|(-A)^{\beta}S_{\alpha}(t)y\|\nonumber\\
&\quad+\|(-A)^{\beta}\int_{0}^{t}P_{\alpha}(t-s)f(s)ds\|.
\end{align}
Put (\ref{change}), (\ref{ca}), (\ref{ap}) into (\ref{varphi}) to get
\begin{align}\label{upper}
\|\varphi(t)-x\|_{\beta}&\leq
C(\alpha,T)\|x\|_{\beta}+M_{0}M'(\alpha,T,y)\nonumber\\
&\quad+M_{0}(M'_{T}Me^{\omega T}T+Me^{\omega T}M_{T}+M_{T})
\end{align}
For $t \in[0,T]$ and $x_{1},x_{2},x_{3}\in N_{\delta}(x)$, we have
\begin{align}\label{following}
&\|(-A)^{\beta-1}\left(-\int_{0}^{t}C_{\alpha}(t-s)(h(s,s,x_{1})+\int_{0}^{s}h_{1}(s,r,x_{2})dr)ds+\int_{0}^{t}h(t,s,x_{3})ds\right)\|\nonumber\\
&\leq M_{0}\{[TMe^{\omega T}(1+M(x,T)+T(1+N(x,T))]+T(1+M(x,T))\}.
\end{align}
From (\ref{upper}) and (\ref{following}), it can be seen that we can choose  $\delta>0$ and $T>0$ such that (\ref{c2}) and (\ref{c3}) hold,
provided that $\delta$ and $T$ satisfy the following inequalities
\begin{align}\label{leq1}
&C(\alpha,T)\|x\|_{\beta}+M_{0}M'(\alpha,T,y)\nonumber\\
&\quad+M_{0}(M'_{T}Me^{\omega T}T+Me^{\omega T}M_{T}+M_{T})<\frac{\delta}{2},
\end{align}
and
\begin{align}\label{leq2}
M_{0}\{[TMe^{\omega T}(1+M(x,T)+T(1+N(x,T))]+T(1+M(x,T))\}<\frac{\delta}{2}.
\end{align}
Let $C:=C([0,T];X_{\beta})$ equipped with the norm $\|\phi\|_{C}=\sup_{t\in [0,T]}\|\phi(t)\|_{\beta}$.  Let $F$ be the closed convex bounded subset of $C([0,T];X_{\beta})$, defined by
\begin{equation*}
F=\{\phi\in C: \|\phi-\varphi\|_{C}\leq \frac{\delta}{2}\}.
\end{equation*}
From $\|\phi(t)-x\|_{\beta}\leq \|\phi-\varphi\|_{C}+\|\varphi(t)-x\|_{\beta}\leq \delta$, it follows that $\phi(t)\in D$ for $\phi(t)\in F$, $t\in [0,T]$. Set the mapping $Q$ on $F$ by
\begin{equation*}
(Q\phi)(t)=\varphi(t)+\int_{0}^{t}P_{\alpha}(t-s)\int_{0}^{s}h(s,r,\phi(r))drds, \ t\in[0,T].
\end{equation*}
Step 1: We show that $Q$ maps $F$ into $F$.
Since
\begin{equation*}
\frac{d}{ds}\int_{0}^{s}h(s,r,\phi(r))dr=\int_{0}^{s}h_{1}(s,r,\phi(r))dr+h(s,s,\phi(s)), \end{equation*}
by (\ref{key}), (\ref{beta}), (\ref{avt}), (\ref{c3}), we have
\begin{align*}
\|(Q\phi)(t)-\varphi(t)\|_{\beta}&=\|(-A)^{\beta}((Q\phi)(t)-\varphi(t))\|\\
&=\|(-A)^{\beta-1}(-\int_{0}^{t}(AP_{\alpha}(t-s)\int_{0}^{s}h(s,r,\phi(r))drds))\|\\
&=\|(-A)^{\beta-1}[-\int_{0}^{t}C_{\alpha}(t-s)(\int_{0}^{s}h_{1}(s,r,\phi(r))dr+h(s,s,\phi(s)))ds\\
&\quad+\int_{0}^{t}h(t,s,\phi(s))ds]\|<\frac{\delta}{2}.
\end{align*}
It is easy to show that $Q\phi: [0,T]\rightarrow X_{\beta}$ is continuous in $t$ on $[0,T]$. We see that $Q$ maps $F$ into $F$.\\
Step 2: We show that $Q$ is continuous. By ($A_{1}$), ($A_{2}$), for every $\varepsilon>0$, there exists some $\delta>0$ such that for $\phi_{1},\phi_{2}\in F$, $\|\phi_{1}-\phi_{2}\|_{C}<\delta$, $s\in [0,T]$,
\begin{align*}
\sup_{r\in [0,T]}\|h(s,r,\phi_{1}(r))-h(s,r,\phi_{2}(r))\|<\varepsilon,\ \sup_{r\in [0,T]}\|h_{1}(s,r,\phi_{1}(r))-h_{1}(s,r,\phi_{2}(r))\|<\varepsilon.
\end{align*}
Then
\begin{align*}
&\|(Q\phi_{1})(t)-(Q\phi_{2})(t)\|_{\beta}\\
&=\|(-A)^{\beta-1}[(-A)\int_{0}^{t}P_{\alpha}(t-s)\int_{0}^{s}h(s,r,\phi_{1}(r))drds\|\\
&\quad-(-A)\int_{0}^{t}P_{\alpha}(t-s)\int_{0}^{s}h(s,r,\phi_{2}(r))drds]\\
&=\|(-A)^{\beta-1}[-\int_{0}^{t}C_{\alpha}(t-s)(\int_{0}^{s}h_{1}(s,r,\phi_{1}(r))dr-\int_{0}^{s}h_{1}(s,r,\phi_{2}(r))dr\\
&\quad+h(s,s,\phi_{1}(s))-h(s,s,\phi_{2}(s)))ds+\int_{0}^{t}(h(t,s,\phi_{1}(s))-\int_{0}^{t}(h(t,s,\phi_{2}(s)))ds]\|\\
&\leq \|(-A)^{\beta-1}\|[\int_{0}^{t}Me^{\omega(t-s)}(\int_{0}^{s}\varepsilon dr+\varepsilon)ds+\int_{0}^{t}\varepsilon ds].
\end{align*}
This implies that $Q$ is continuous.\\
Step 3. We show that the set $\{Q\phi:\phi\in F\}$ is equicontinuous. For $\phi\in F$, $0\leq t\leq t'\leq T$, we have
\begin{align*}
\|(Q\phi)(t)-(Q\phi)(t')\|_{\beta}&\leq \|(C_{\alpha}(t)-C_{\alpha}(t'))(-A)^{\beta}x\|+\|A(P_{\alpha}(t)-P_{\alpha}(t'))(-A)^{\beta-1}y\|\\
&\quad+\|(-A)^{\beta-1}[\int_{0}^{t}(C_{\alpha}(t-s)(h(s,s,\phi(s)+\int_{0}^{s}h_{1}(s,r,\varphi(r))dr)ds\\
&\quad-\int_{0}^{t'}C_{\alpha}(t'-s)(h(s,s,\phi(s)+\int_{0}^{s}h_{1}(s,r,\varphi(r))dr)ds]\|\\
&\quad+\|(-A)^{\beta-1}(\int_{0}^{t}h(t,s,\phi(s)ds-\int_{0}^{t'}h(t',s,\phi(s))ds)\|\\
&\quad+\|(-A)^{\beta-1}(\int_{0}^{t}C_{\alpha}(t-s)f'(s)ds-\int_{0}^{t'}C_{\alpha}(t'-s)f'(s)ds)\|\\
&\quad+\|(-A)^{\beta-1}(C_{\alpha}(t)-C_{\alpha}(t'))f(0)\|+\|(-A)^{\beta-1}(f(t)-f(t'))\|.
\end{align*}
Since $C_{\alpha}(t)$, $P_{\alpha}(t)$ is strongly continuous, it follows that
\begin{align*}
\|(C_{\alpha}(t)-C_{\alpha}(t'))(-A)^{\beta}x\|+\|A(P_{\alpha}(t)-P_{\alpha}(t'))(-A)^{\beta-1}y\|\rightarrow 0
\end{align*}
as $|t-t'|\rightarrow 0$, and
\begin{align*}
\|(-A)^{\beta-1}(C_{\alpha}(t)-C_{\alpha}(t'))f(0)\|+\|(-A)^{\beta-1}(f(t)-f(t'))\|\rightarrow 0
\end{align*}
as $|t-t'|\rightarrow 0$.\\
By Lemma 2.1 in \cite{Travis}, since $A^{-1}$ is compact, then for $0<\beta<1$, $(-A)^{\beta-1}$ is compact.
The compactness of $(-A)^{\beta-1}$, the strongly continuity of $C_{\alpha}(t), P_{\alpha}(t)$, together with (\ref{c1}) imply that
\begin{align*}
&\|(-A)^{\beta-1}[\int_{0}^{t}(C_{\alpha}(t-s)(h(s,s,\phi(s)+\int_{0}^{s}h_{1}(s,r,\varphi(r))dr)ds\\
&\quad-\int_{0}^{t'}C_{\alpha}(t'-s)(h(s,s,\phi(s))+\int_{0}^{s}h_{1}(s,r,\varphi(r))dr)ds]\|\\
&\leq \|\int_{0}^{t}(C_{\alpha}(t-s)-C_{\alpha}(t'-s))(-A)^{\beta-1}(h(s,s,\phi(s)+\int_{0}^{s}h_{1}(s,r,\varphi(r))dr)ds\|\\
&\quad+\|(-A)^{\beta-1}\|\|\int_{t}^{t'}C_{\alpha}(t'-s)(h(s,s,\phi(s))+\int_{0}^{s}h_{1}(s,r,\varphi(r))dr)ds\|\rightarrow 0
\end{align*}
as $|t-t'|\rightarrow 0$.\\
On the other hand, by (\ref{c1}),
\begin{align*}
&\|(-A)^{\beta-1}(\int_{0}^{t}C_{\alpha}(t-s)f'(s)ds-\int_{0}^{t'}C_{\alpha}(t'-s)f'(s)ds)\|\\
&\leq \|(-A)^{\beta-1}\|(\|\int_{0}^{t}\int_{t}^{t'}h_{1}(r,s,\phi(s))drds\|+\|\int_{t}^{t'}h(t',s,\phi(s))ds\|)\rightarrow 0
\end{align*}
as $|t-t'|\rightarrow 0$, and
\begin{align*}
&\|(-A)^{\beta-1}(\int_{0}^{t}C_{\alpha}(t-s)f'(s)ds-\int_{0}^{t'}C_{\alpha}(t'-s)f'(s)ds)\|\\
&\leq \|(-A)^{\beta-1}\|(\|\int_{0}^{t}(C_{\alpha}(t-s)-C(t'-s))f'(s)ds\|+\|\int_{t}^{t'}C_{\alpha}(t'-s)f'(s)ds\|)\rightarrow 0
\end{align*}
as $|t-t'|\rightarrow 0$.\\
Therefore, $\{Q\phi:\phi\in F\}$ is equicontinuous.\\
Step 4. We show that for any given $t\in[0,T]$, the set $\{Q\phi:\phi\in F\}$ is precompact in $X_{\beta}$. Since $A^{-1}$ is compact, then for $\gamma\in(\beta,1]$, $(-A)^{-\gamma}:X\rightarrow X_{\beta}$ is compact, we only need to prove that
$\{(-A)^{\gamma}((Q\phi)(t)-\varphi(t)):\phi\in F\}$ is bounded $\gamma\in(\alpha,1]$.
In fact, we have
\begin{align*}
\|(-A)^{\gamma}(Q\phi-\varphi)(t)\|&\leq  \|(-A)^{\gamma-1}\int_{0}^{t}C_{\alpha}(t-s)(h(s,s,\phi(s))+\int_{0}^{s}h_{1}(s,r,\varphi(r))dr)ds\\
&\quad+(-A)^{\gamma-1}\int_{0}^{t}h(t,s,\phi(s))\|.
\end{align*}
From (\ref{c1}), the boundedness is obtained. Therefore, it follows from the Arzela-Ascoli theorem, $Q$ is compact. By Schauder fixed point theorem, $Q$ has a fixed point in $F$, which is a solution of (\ref{ut}). If $x\in D(A)$, $y\in E$, then by Lemma 3.4, the solution of (\ref{ut}) is a solution of (\ref{Caputo}). \ \ \ $\Box$

\section{An example}
Consider the fractional  semilinear Volterra integrodifferential equation of order $\alpha\in(1,2]$
\begin{equation}\label{Volterra}
 \left\{\begin{aligned} &^{C}D_{t}^{\alpha}z(t,x)=\Delta z(t,x)+\int_{0}^{t}\rho(t,s,z(s,x))ds+\theta(t,x), \ t\in R_{+},\ x\in(0,\pi),\\
&z(t,0)=z(t,\pi), \ t\in R_{+},\\
&z(0,x)=\sigma(x),\ z_{t}(0,x)=\mu(x), \  x\in(0,\pi),
\end{aligned}\right.
\end{equation}
where $^{C}D_{t}^{\alpha}$ is the $\alpha$-order Caputo fractional derivative operator. Let $X=L^{2}[0,\pi]$ and define $A:X\rightarrow X$ by $Aw=w''$ with the domain
$D(A)=\{w\in X: w,w' \ \mbox{are absolutely continuous}, \ w''\in X, \ w(0)=w(\pi)=0\}$.
Thus
\begin{equation*}
Aw=-\sum_{n=1}^{\infty}n^{2}(w,w_{n})w_{n},\ w\in D(A),
\end{equation*}
where $w_{n}(s)=\sqrt{\pi/2}\sin ns$, $n=1,2,\cdots$, is the orthonormal set of eigenvalues of $A$. It is easy to see that $A$ is the infinitesimal generator of a strongly continuous cosine family $C(t), \ t\in R$ on $X$ given by $C(t)w=\sum_{n=1}^{\infty}\cos nt(w,w_{n})w_{n}$, $w\in X$.
From the subordinate principle (see Theorem 3.1 in \cite{bazh}), it follows that $A$ is the infinitesimal generator of $\alpha$-order cosine family $C_{\alpha}(t)$ such that $C_{\alpha}(0)=I$, and
\begin{equation*}
C_{\alpha}(t)=\int_{0}^{\infty}\varphi_{t,\alpha/2}(s)C(s)ds, \quad t>0,
\end{equation*}
where $\varphi_{t,\alpha/2}(s)=t^{-\alpha/2}\phi_{\alpha/2}(st^{-\alpha/2})$, and
\begin{equation*}
\phi_{\gamma}(z)=\sum_{n=0}^{\infty}\frac{(-z)^{n}}{n!\Gamma(-\gamma n+1-\gamma)}, \quad 0<\gamma<1.
\end{equation*}
If we take $\beta=\frac{1}{2}$, then
\begin{equation*}
(-A)^{1/2}w=\sum_{n=1}^{\infty}n(w,w_{n})w_{n},\ w\in D((-A)^{1/2}).
\end{equation*}
The operator $(-A)^{-1/2}$ is given by
\begin{equation*}
(-A)^{-1/2}w=\sum_{n=1}^{\infty}(1/n)(w,w_{n})w_{n},\ w\in X.
\end{equation*}
It is easy to show that $(-A)^{-1/2}$ is compact. By Lemma 2.1 in \cite{Travis}, $A^{-1}$ is compact.
Let $\rho: R_{+}\times R_{+}\times R_{+} \rightarrow R$ be continuous and continuously differentiable with respect to its first variable. Let
$\theta: R_{+}\times R_{+}\rightarrow R$ be continuous and continuously differentiable with respect to its first variable. Let $h: R_{+}\times R_{+} \times X_{1/2}\rightarrow X$ be defined by $(h(t,s,w))(x)=\rho(t,s,w(x))$, $w\in X_{1/2}$, $x\in [0,\pi]$, and let $f:R_{+}\rightarrow X$ be defined by
$(f(t))(x)=\theta(t,x)$, $x\in [0,\pi]$. Then we can rewrite (\ref{Volterra}) as (\ref{Caputo}). If $w\in D((-A)^{1/2})$, then $w$ is absolutely continuous, $w'\in X$, $w(0)=w(\pi)=0$, and $\|w\|_{1/2}=\|w'\|$ (see Chapter 6 in \cite{Milan}). Let $t_{1},s_{1}\in [0,T]$, $w_{1}\in X_{1/2}$. For every $\varepsilon>0$,  there exists a $\delta>0$ such that if $t,s\in [0,T]$, $x\in[0,\pi]$, $v\in R$, and $|t_{1}-t|<\delta$,
$|s_{1}-s|<\delta$, $|w_{1}(x)-v|<\delta$, then $|\rho(t_{1},s_{1},w_{1}(x))-\rho(t,s,v)|<\varepsilon$. Let $w\in X_{1/2}$, and $\|w_{1}-w\|_{1/2}<\delta/\sqrt{\pi}$. Then $|w_{1}(x)-w(x)|\leq |\int_{0}^{x}(w_{1}'(r)-w'(r))dr|\leq \int_{0}^{x}|w_{1}'(r)-w'(r)|dr\leq\sqrt{\pi}\|w_{1}'-w'\|=\sqrt{\pi}\|w_{1}-w\|_{1/2}$.
Hence, for  $|t_{1}-t|<\delta$,
$|s_{1}-s|<\delta$, $|w_{1}(x)-v|<\delta$, we have $\|h(t_{1},s_{1},w)-h(t,s,w)\|=\int_{0}^{\pi}|\rho(t_{1},s_{1},w_{1}(x)-\rho(t,s,w(x))|^{2}dx\leq\pi\varepsilon^{2}$.
Therefore $h$ is continuous. By similar method, the conditions $(A_{2})$ and $(A_{3})$ are satisfied. By Theorem 3.6, the integrodifferential equation (\ref{Volterra}) has a local solution.


\begin{thebibliography}{99}
\bibitem{Pruss} J. Pr$\ddot{u}$ss,  Bounded solutions of Volterra equations, SIAM J. Math. Anal. 19 (1988) 133-149.
\bibitem{Jan} J. Pr$\ddot{u}$ss, Positivity  and regularity
of hyperbolic Volterra equations in Banach spaces, Math. Ann. 279 (1987) 317-344.
\bibitem{J} J. Pr$\ddot{u}$ss, On linear Volterra equations of parabolic type in Banach spaces, Trans. Amer. Math. Soc. 301 (1987) 691-721.
\bibitem{Grimmer} R. C. Grimmer, Resolvent operators for integral equations in a Banach space, Trans. Amer. Math. Soc. 273 (1982) 333-349.
\bibitem{WEF} W. Fitzgibbon, `Semilinear integrodifferential equations in Banach space', Nonlinear Anal. 4 (1980) 745-760.
\bibitem{Valentin} V. Keyantuo and C. Lizama, H$\ddot{o}$lder continuous solutions for integro-differential equations and maximal regularity, J. Differential. Equ. 230 (2006) 634-660.
\bibitem{Londen} S. O. Londen, An existence result on a volterra equation in a Banach space, Trans. Amer. Math. Soc. 235 (1978) 285-304.
\bibitem{MM} E. Mainini and  G. Mola, Exponential and polynomial decay for first order
linear Volterra equations, Quart. Appl. Math. Volume LXVII (2009) 93-111.
\bibitem{bazh} E. Bazhlekova, Fractional Evolution Equations in Banach
Spaces, PhD Thesis, Eindhoven University of
Technology, 2001.
\bibitem{Pazy} A. Pazy, Semigroups of Linear Operators and Applications to
Partial Differential Equations (Springer-Verlag, New York, 1983).
\bibitem{Pod} I. Podlubny, Fractional Differential Equations (
Academic Press, New York, 1999).
\bibitem{Chenchuang} C. Chen and M. Li, On
fractional resolvent operator functions. Semigroup Forum, 80 (2010) 121-142.
\bibitem{Peng} J. Peng and K. Li,  A novel characteristic of
solution operator for the fractional abstract Cauchy problem,  J.
Math. Anal. Appl. 385 (2012) 786-796.
\bibitem{Miao} M. Li, C. Chen and  F.B. Li, On fractional powers of generators of
fractional resolvent families, J. Funct. Anal. 259 (2010) 2702-2726.
\bibitem{Travis} C. C. Travis and G. F. Webb, An abstract second order semilinear volterra integrodifferential equation, SIAM J. Math. Anal. 10 (1979) 412-424.
\bibitem{Milan} M. Miklav$\hat{c}$i$\hat{c}$, Applied Functional Analysis and Partial Differential Equations (World Scientific, Singapore, 1998).
\bibitem{Main} F. Mainardi, Fractional Calculus and Waves in Linear Viscoelasticity: An Introduction to Mathematical Models (Imperial College Press, London 2010).
\bibitem{Lang} C. Lang and J. Chang, Local existence for nonlinear Volterra integrodifferential equations with infinite delay, Nonlinear Anal. 68 (2008) 2943-2956.
\bibitem{Adel} A. Jawahdou,  Mild solutions of functional semilinear evolution Volterra
integrodifferential equations on an unbounded interval, Nonlinear Anal. 74 (2011) 7325-7332.
\bibitem{EH} E. Hernandez, $\mathcal{C}^{\alpha}$-classical solutions for abstract non-autonomous integro-differential equations, Proc. Amer. Math. Soc 139 (2011) 4307-4318.
\bibitem{Keyantuo} V. Keyantuo and C. Lizama, On a connection between powers of operators and fractiona Cauchy problems, J. Evol. Equ. 12 (2012) 245-265.
\bibitem{Engler} H. Engler, Weak  solutions of a  class of quasilinear
hyperbolic  integro-differential equations describing  viscoelastic materials, Arch. Rational Mech. Anal. 113 (1991) 1-38.
\bibitem{R} R. Metzler and J. Klafter, The random walk's guide to anomalous diffusion: a fractional dynamics approach, Phys. Rep. 339 (2000) 1-77.
\end{thebibliography}
\end{document}